\documentclass[12pt,a4paper]{article}
\usepackage{lineno,hyperref}
\usepackage[latin5]{inputenc}
\usepackage{amsmath}
\usepackage{amssymb}
\usepackage{graphicx}
\usepackage{subfigure}
\usepackage{caption}
\usepackage{float}
\usepackage{mathtools}
\usepackage{color}
\usepackage{mathrsfs}
\usepackage[]{geometry}
\usepackage{geometry}
\usepackage{subfigure}

\title{Numerical Solutions of the Gardner Equation by Extended Form of the Cubic B-splines}
\author{Ozlem Ersoy Hepson$^{a,}$\thanks{Corresponding Author: ozersoy@ogu.edu.tr}, Alper Korkmaz$^{b,}$, Idiris Dag$^c$\\
{\scriptsize $^{a}$ Eskisehir Osmangazi University, Mathematics \& Computer Department, Eskisehir, Turkey.}\\
   {\scriptsize          $^{b}$ Çankırı Karatekin University, Department of Mathematics, Çankırı, Turkey.}\\
   {\scriptsize          $^{c}$ Eskisehir Osmangazi University, Computer Engineering Department, Eskisehir, Turkey.}
   } 
\begin{document}
\maketitle
\begin{abstract}
The extended definition of the polynomial B-splines may give a chance to improve the results obtained by the classical cubic polynomial B-splines. Determination of the optimum value of the extension parameter can be achieved by scanning some intervals containing zero. This study aims to solve some initial boundary value problems constructed for the Gardner equation by the extended cubic B-spline collocation method. The test problems are derived from some analytical studies to validate the efficiency and accuracy of the suggested method. The conservation laws are also determined to observe them remain constant as expected in theoretical aspect. The stability of the proposed method is investigated by the Von Neumann analysis.
\end{abstract}
\section{Introduction}

Consider the Gardner equation (or combined KdV-mKdV) of the form 
\begin{equation}
u_{t}+\mu_1 uu_{x}+\mu_2 u^2u_{x}+\mu_3 u_{xxx}=0 \label{gardner}
\end{equation}
where $u=u(x,t)$ and $\mu_1$, $\mu_2$ and $\mu_3$ are constants. The Gardner equation has two nonlinear terms in the quadratic and cubic forms and the dissipative term is of third order. The Gardner equation is an integrable system and Miura transformation connects it to the KdV equation\cite{demler1}. 

The Gardner equation is a useful model to understand the propagation of negative ion acoustic plasma waves\cite{ruderman1}. The equation can be derived from the system of plasma motion equations in one dimension with arbitrarily charged cold ions and inertia neglected isothermal electrons.    

The Gardner equation can also be a good description of internal waves with large amplitudes\cite{kamc1}. The modulation system in the Riemann invariant form and classified solutions for large $t$ values subject to some particular initial conditions are studied deeply to explain the undular bore formation. The structure of the nonlinear terms and the positions of the initial step are significant for the solution classification. The sign of $\mu_2$ has also critical role for the structures of the solutions. The patterns of the solutions can be of the forms bright or dark cnoidal or trigonometric bores, kinks, rarefaction waves or combinations of them.  

Some particular forms of the internal ocean rogue waves in the coasts occur in the Gardner equation when the $\mu_2$ and $\mu_3$ have the same sign which enables modulational instability\cite{grimshaw1}. The Gardner equation can also be a useful model to study stratified fluid transcritical flow passing a topographic obstacle
when a forcing term is added to the equation\cite{kamc2}.

Interacting two-soliton type solution is derived by Darboux transformations in Slyunyaev and Pelinovskii's study\cite{sul1}. The consistent tanh method is also capable of generating interacting solutions for the Gardner equation\cite{hu1}. Some interacting of two wave solutions such as soliton-cnoidal wave or soliton-periodic wave are presented in this study. The consistent Riccati expansion is another method that is capable to obtain soliton-cnoidal wave interaction type solutions\cite{feng1}.   

Dynamics of solitons, in the adiabatic parameter case, of the perturbed Gardner equation is studied in details by Biswas and Zerrad\cite{biswas1}. The motion integrals are also determined for the perturbed form of the equation. The effects of perturbation to a single soliton type solution are discussed by considering the $\tanh$ type initial data\cite{jia1}.  

Exact solutions to the Gardner equation are set up by using various methods. some solutions containing $\tanh$ and $\coth$ functions are proposed by the extended form of the $\tanh$ method\cite{bekir1}. Solitary wave and periodic solutions are constructed by aid of the projective Riccati equations\cite{fu1}. These solutions have various terms including trigonometric or hyperbolic functions in rational forms. $G'/G$ is another expansion method to determine the exact solutions of the Gardner equation. Some solitary wave, periodic, exponential, rational and complex-type traveling wave solutions are found by this method\cite{lu1}. Some other exact solutions in terms of trigonometric functions \cite{akbar1,naher1,jawad1}, hyperbolic functions\cite{akbar1,naher1,taghizade1}, kink solutions\cite{wazwaz1} are determined by using different expansion or ansatz methods. 

Numerical solutions to the Gardner equations are also concerned in various studies. The conservative finite difference schemes are developed to determine propagation of one soliton and collusion of two soliton solutions numerically\cite{nishiyama1}. Restrictive Taylor's technique has also been implemented to simulate the propagation of some solutions numerically\cite{rageh1}.

In the present study, we develop an extended cubic B-spline collocation method to the solutions of the Gardner equation. Having only two continuous derivatives of the cubic B-splines force us to reduce the order of the third order derivative term. Setting $v=u_x$ reduces the order of the Gardner equation (\ref{gardner}) in the resultant coupled system of nonlinear PDEs of the form
\begin{equation}
\begin{aligned}
u_{t}+\left(\mu_1 u+\mu_2 u^{2}\right)u_x+\mu_3 v_{xx}&=0 \\ 
v-u_{x}&=0
\end{aligned}
\label{2}
\end{equation}
In order to complete the mathematical notation of the initial-boundary value problems (IVPs), we use the initial data
\begin{equation}
\begin{aligned}
u(x,0)&=f(x) \\
v(x,0)&=f_x(x)
\end{aligned}
\end{equation}
and the homogeneous Neumann conditions
\begin{equation}
\begin{aligned}
u_x(a,t)&=0,\, u_x(b,t)&=0, \\
u_{xx}(a,t)&=0,\, u_{xx}(b,t)&=0, \\
v_x(a,t)&=0,\, v_x(b,t)&=0, \\
v_{xx}(a,t)&=0,\, v_{xx}(b,t)&=0 
\end{aligned}
\end{equation}
at both ends of the problem interval $[a,b]$. 
\section{Numerical Approximate}
Consider the equal grid distribution 
\begin{equation*}
\pi :a=x_{0}<x_{1}<\ldots <x_{N}=b
\end{equation*}%
of the finite interval $[a,b]$ where $h=(b-a)/N$ is the equal mesh size. An extended B-spline function is defined as
\begin{equation}
E_{j}(x)=\frac{1}{24h^{4}}\left \{ 
\begin{array}{ll}
4h(1-\lambda )(x-x_{j-2})^{3}+3\lambda (x-x_{j-2})^{4}, & \left[
x_{j-2},x_{j-1}\right] , \\ 
\begin{array}{l}
(4-\lambda )h^{4}+12h^{3}(x-x_{j-1})+6h^{2}(2+\lambda )(x-x_{j-1})^{2} \\ 
-12h(x-x_{j-1})^{3}-3\lambda (x-x_{j-1})^{4}%
\end{array}
& \left[ x_{j-1},x_{i}\right] , \\ 
\begin{array}{l}
(4-\lambda )h^{4}-12h^{3}(x-x_{j+1})+6h^{2}(2+\lambda )(x-x_{j+1})^{2} \\ 
+12h(x-x_{j+1})^{3}-3\lambda (x-x_{j+1})^{4}%
\end{array}
& \left[ x_{i},x_{j+1}\right] , \\ 
4h(\lambda -1)(x-x_{j+2})^{3}+3\lambda (x-x_{j+2})^{4}, & \left[
x_{j+1},x_{j+2}\right] , \\ 
0 & \text{otherwise.}%
\end{array}%
\right.  \label{e1}
\end{equation}
with the extension parameter $\lambda$. The set of extended B-spline functions $\{E_{-1}(x),E_{0}(x),\dots ,E_{N+1}(x)\}$ constitutes a basis function set for the functions defined in this interval. The approximate solutions $U(x,t)$ and $V(x,t)$ to $u(x,t)$ and $v(x,t)$, respectively, can be written in terms of the extended B-splines as
\begin{equation}
U(x,t)=\sum_{j=-1}^{N+1}\delta _{j}E_{j}(x),\text{ }V(x,t)=%
\sum_{j=-1}^{N+1}\phi _{j}E_{j}(x)  \label{e2}
\end{equation}
where $\delta _{j}$ and $\phi _{j}$ are time dependent parameters. These parameters are determined after implementation of the collocation method and complementary data. Since each extended B-spline has lowest two derivatives, the nodal derivative values of both $U$ and $V$ can be summarized in terms of extended B-splines as
\begin{equation}
\begin{tabular}{l}
$U_{j}=U(x_{j},t)=\dfrac{4-\lambda }{24}\delta _{j-1}+\dfrac{8+\lambda }{12}%
\delta _{j}+\dfrac{4-\lambda }{24}\delta _{j+1},$ \\ 
\\ 
$U_{j}^{\prime }=U^{\prime }(x_{j},t)=\dfrac{-1}{2h}\left( \delta
_{j-1}-\delta _{j+1}\right) $ \\ 
\\ 
$U_{j}^{\prime \prime }=U^{\prime \prime }(x_{j},t)=\dfrac{2+\lambda }{2h^{2}%
}\left( \delta _{j-1}-2\delta _{j}+\delta _{j+1}\right) $%
\end{tabular}%
\begin{tabular}{l}
$V_{j}=V(x_{j},t)=\dfrac{4-\lambda }{24}\phi _{j-1}+\dfrac{8+\lambda }{12}%
\phi _{j}+\dfrac{4-\lambda }{24}\phi _{j+1},$ \\ 
\\ 
$V_{j}^{\prime }=V^{\prime }(x_{j},t)=\dfrac{-1}{2h}\left( \phi _{j-1}-\phi
_{j+1}\right) $ \\ 
\\ 
$V_{j}^{\prime \prime }=V^{\prime \prime }(x_{j},t)=\dfrac{2+\lambda }{2h^{2}%
}\left( \phi _{j-1}-2\phi _{j}+\phi _{j+1}\right) $%
\end{tabular}
\label{e3}
\end{equation}

The time integration of the space-splitted system (\ref{2}) is performed by
the Crank-Nicolson method as%
\begin{equation}
\begin{array}{r}
\dfrac{U^{n+1}-U^{n}}{\Delta t}+\mu _{1}\dfrac{(UU_x)^{n+1}+(UU_x)^{n}}{2}+\mu
_{2}\dfrac{(U^{2}U_x)^{n+1}+(U^{2}U_x)^{n}}{2}+\mu _{3}\dfrac{%
V_{xx}^{n+1}+V_{xx}^{n}}{2}=0 \\ 
\\ 
\dfrac{V^{n+1}+V^{n}}{2}-\dfrac{U_{x}^{n+1}+U_{x}^{n}}{2}=0%
\end{array}
\label{10}
\end{equation}%
where the superscript $^p$ represents the solution at the $p$th
time level with equal time step size $\Delta t$ satisfying $t^{n+1}=t^{n}+\Delta t$.

The nonlinear terms $(UU_x)^{n+1}$ and $\left( U^{2}U_x\right) ^{n+1}$ in Eq. (\ref{10}) are converted to linear forms by using 
\begin{equation*}
(UU_x)^{n+1}=U^{n+1}U_x^{n}+U^{n}U_x^{n+1}-U^{n}U_x^{n}
\end{equation*}%
and%
\begin{equation*}
(U^{2}U_x)^{n+1}=2U^{n+1}U^{n}U_x^{n}+(U^{n})^{2}U_x^{n+1}-2(U^{n})^{2}U_x^{n}
\end{equation*}
defined in \cite{rubin}. The resultant linear system is discretized in time by using Crank-Nicolson method as\scriptsize{
\begin{equation}
\begin{aligned}
&\left[ \left( \frac{2}{\Delta t}+\mu _{1}L+2\mu _{2}KL\right) \alpha _{1}+\left( \mu _{1}K+\mu _{2}K^{2}\right)
\beta _{1}
\right] \delta _{j-1}^{n+1}+\left[ \mu _{3}\gamma _{1}\right] \phi _{j-1}^{n+1}+\left[ \left( \frac{%
2}{\Delta t}+\mu _{1}L+2\mu _{2}KL\right) \alpha _{2}\right] \delta
_{j}^{n+1}\\
&+\left[ \mu_{3}\gamma _{2}\right] \phi _{j}^{n+1} \left[ \left( \frac{2}{\Delta t}+\mu _{1}L+2\mu _{2}KL\right) \alpha _{1}-\left( \mu _{1}K+\mu _{2}K^{2}\right)
\beta _{1}\right] \delta _{j+1}^{n+1}+\left[ \mu _{3}\gamma _{1}\right] \phi _{j+1}^{n+1}\\
&=\left[ \left( \frac{2}{\Delta t}+\mu _{2}KL\right) \alpha _{1}\right]
\delta _{j-1}^{n}-\mu _{3}\gamma _{1}\phi _{j-1}^{n}+\left[ \left( \frac{2}{%
\Delta t}+\mu _{2}KL\right) \alpha _{2}\right] \delta _{j}^{n}-\mu
_{3}\gamma _{2}\phi _{j}^{n}+\left[ \left( \frac{2}{\Delta t}+\mu
_{2}KL\right) \alpha _{1}\right] \delta _{j+1}^{n}&-\mu _{3}\gamma _{1}\phi
_{j+1}^{n} \\
&-\beta _{1}\delta _{j-1}^{n+1}+\alpha _{1}\phi _{j-1}^{n+1}+\alpha _{2}\phi
_{j}^{n+1}+\beta _{1}\delta _{j+1}^{n+1}+\alpha _{1}\phi _{j+1}^{n+1}
=\beta _{1}\delta _{j-1}^{n}-\alpha _{1}\phi _{j-1}^{n}-\alpha _{2}\phi
_{j}^{n}-\beta _{1}\delta _{j+1}^{n}-\alpha _{1}\phi _{j+1}^{n} \\
&m=0,...,N,\quad n=0,1...,
\end{aligned}\label{12}
\end{equation}}
\normalsize
where%
\begin{equation*}
\begin{array}{l}
K=\alpha _{1}\delta _{j-1}^{n}+\alpha _{2}\delta _{j}^{n}+\alpha _{1}\delta
_{j+1}^{n} \\ 
L=\alpha _{1}\phi _{j-1}^{n}+\alpha _{2}\phi _{j}^{n}+\alpha _{1}\phi
_{j+1}^{n}%
\end{array}%
\end{equation*}%
\begin{eqnarray*}
\alpha _{1} &=&\dfrac{4-\lambda }{24},\text{ }\alpha _{2}=\dfrac{8+\lambda }{%
12} \\
\beta _{1} &=&-\dfrac{1}{2h},\text{ }\gamma _{1}=\dfrac{2+\lambda }{2h^{2}},%
\text{ }\gamma _{2}=-\dfrac{4+2\lambda }{2h^{2}}
\end{eqnarray*}%
This system can be written in the matrix notation as 
\begin{equation}
\mathbf{Ax}^{n+1}=\mathbf{Bx}^{n}  \label{13}
\end{equation}%
where%
\begin{equation*}
\mathbf{A=}%
\begin{bmatrix}
\nu _{m1} & \nu _{m2} & \nu _{m3} & \nu _{m4} & \nu _{m5} & \nu _{m2} &  & 
&  &  \\ 
-\beta _{1} & \alpha _{1} & 0 & \alpha _{2} & \beta _{1} & \alpha _{1} &  & 
&  &  \\ 
&  & \nu _{m1} & \nu _{m2} & \nu _{m3} & \nu _{m4} & \nu _{m5} & \nu _{m2} & 
&  \\ 
&  & -\beta _{1} & \alpha _{1} & 0 & \alpha _{2} & \beta _{1} & \alpha _{1}
&  &  \\ 
&  &  & \ddots  & \ddots  & \ddots  & \ddots  & \ddots  & \ddots  &  \\ 
&  &  &  & \nu _{m1} & \nu _{m2} & \nu _{m3} & \nu _{m4} & \nu _{m5} & \nu
_{m2} \\ 
&  &  &  & -\beta _{1} & \alpha _{1} & 0 & \alpha _{2} & \beta _{1} & \alpha
_{1}%
\end{bmatrix}%
\end{equation*}

\begin{equation*}
\mathbf{B=}%
\begin{bmatrix}
\nu _{m6} & -\nu _{m2} & \nu _{m7} & -\nu _{m4} & \nu _{m6} & -\nu _{m2} &  & 
&  &  \\ 
\beta _{1} & -\alpha _{1} & 0 & -\alpha _{2} & -\beta _{1} & -\alpha _{1} & 
&  &  &  \\ 
&  & \nu _{m6} & -\nu _{m2} & \nu _{m7} & -\nu _{m4} & \nu _{m6} & -\nu _{m2} & 
&  \\ 
&  & \beta _{1} & -\alpha _{1} & 0 & -\alpha _{2} & -\beta _{1} & -\alpha
_{1} &  &  \\ 
&  &  & \ddots & \ddots & \ddots & \ddots & \ddots & \ddots &  \\ 
&  &  &  & \nu _{m6} & -\nu _{m2} & \nu _{m7} & -\nu _{m4} & \nu _{m6} & -\nu _{m2} \\ 
&  &  &  & \beta _{1} & -\alpha _{1} & 0 & -\alpha _{2} & -\beta _{1} & 
-\alpha _{1}%
\end{bmatrix}%
\end{equation*}%
and%
\begin{equation*}
\begin{array}{ll}
\nu _{m1}=\left( \frac{2}{\Delta t}+\mu _{1}L+2\mu _{2}KL\right) \alpha _{1}+\left( \mu _{1}K+\mu _{2}K^{2}\right)
\beta _{1}
& \nu _{m6}=\left( \frac{2}{\Delta t}+\mu _{2}KL\right) \alpha _{1} \\ 
\nu _{m2}=\mu _{3}\gamma_{1} &  \nu _{m5}=\left( \frac{2}{\Delta t}+\mu _{1}L+2\mu _{2}KL\right) \alpha _{2}-\left( \mu _{1}K+\mu _{2}K^{2}\right)
\beta _{1}\\ 
\nu _{m3}=\left( \frac{2}{\Delta t}+\mu _{1}L+2\mu _{2}KL\right) \alpha _{2}-\left( \mu _{1}K+\mu _{2}K^{2}\right)
\beta _{1}
& \nu _{m7}=\left( \frac{2}{\Delta t}+\mu _{2}KL\right) \alpha _{2} \\ 
\nu _{m4}=\mu _{3}\gamma
_{2} & 
\end{array}%
\end{equation*}

The system (\ref{13}) consists of $2N+2$ linear equations and $2N+6$ unknowns 
\begin{equation*}
\mathbf{x}^{n+1}=(\delta _{-1}^{n+1},\phi _{-1}^{n+1},\delta _{0}^{n+1},\phi
_{0}^{n+1},\ldots ,\delta _{n+1}^{n+1},\phi _{n+1}^{n+1},).
\end{equation*}
The unique solution of this system requires additional four constraints. The boundary data $U_{x}(a,t)=0,$ $V_{x}(a,t)=0$ and $U_{x}(b,t)=0,$ $V_{x}(b,t)=0$  can be written in terms of parameters as 
the following equations:
\begin{equation*}
\begin{array}{l}
\delta _{-1}=\delta _{1} \\ 
\phi _{-1}=\phi _{1} \\ 
\delta _{N-1}=\delta _{N+1} \\ 
\phi _{N-1}=\phi _{N+1}%
\end{array}%
\end{equation*}

The parameters $\delta _{-1},\phi _{-1},\delta _{N+1},\phi
_{N+1}$ in Eq.(\ref{13}) are eliminated from the system by using the boundary data to determine a solvable system. In order to initialize the iteration algorithm, initial parameters $\delta _{j}^{0},\phi _{j}^{0},$ $j=-1,\ldots ,N+1$ are determined by using the data obtained from the complementary data as 
\begin{equation*}
\begin{array}{l}
U_{x}(a,0)=0=\delta _{-1}^{0}-\delta _{1}^{0}, \\ 
U(x_{j},0)=\alpha _{1}\delta _{j-1}^{0}+\alpha _{2}\delta _{j}^{0}+\alpha
_{1}\delta _{j+1}^{0}=u(x_{j},0),j=1,...,N-1 \\ 
U_{x}(b,0)=0=\delta _{N-1}^{0}-\delta _{N+1}^{0}, \\ 
V_{x}(a,0)=0=\phi _{-1}^{0}-\phi _{1}^{0} \\ 
V(x_{j},0)=\alpha _{1}\phi _{j-1}^{0}+\alpha _{2}\phi _{j}^{0}+\alpha
_{1}\phi _{j+1}^{0}=v(x_{j},0),j=1,...,N-1 \\ 
V_{x}(a,0)=\phi _{N-1}^{0}-\phi _{N+1}^{0}%
\end{array}%
\end{equation*}

\section{Stability Analysis}
The stability of the method is investigated by performing the Von-Neumann analysis where 
\begin{eqnarray}
\delta_{j}^{n} &=&A_1\xi ^{n}\exp (ij\varphi )  \label{k} \\
\phi _{j}^{n} &=&A_2\xi ^{n}\exp (ij\varphi )  \notag
\end{eqnarray}%
\begin{equation*}
\rho=\frac{\xi ^{n+1}}{\xi ^{n}}
\end{equation*}%
Here, $A_1$ and $A_2$ represent the harmonics amplitude. $k$ is the
mode number, $\rho$ is the amplification factor and $\varphi =kh$. The term $U+U^2$ is assumed as locally constant and replaced $\varepsilon$. Substituting \ref{k} into the system 
\begin{eqnarray}
&&a_{1}\delta_{j-1}^{n+1}+a_{2}\delta_{j}^{n+1}+a_{1}\delta_{j+1}^{n+1}+\frac{\lambda
k\varepsilon }{2}(a_{3}\delta_{j-1}^{n+1}-a_{3}\delta_{j+1}^{n+1})+\frac{k\mu_3 }{2}%
(a_{4}\phi _{j-1}^{n+1}+a_{5}\phi _{j}^{n+1}+a_{4}\phi _{j+1}^{n+1})
\label{k1} \\
&=&a_{1}\delta_{j-1}^{n}+a_{2}\delta_{j}^{n}+a_{1}\delta_{j+1}^{n}-\frac{\lambda
k\varepsilon }{2}(a_{3}\delta_{j-1}^{n}-a_{3}\delta_{j+1}^{n})-\frac{k\mu_3 }{2}%
(a_{4}\phi _{j-1}^{n}+a_{5}\phi _{j}^{n}+a_{4}\phi _{j+1}^{n})  \notag
\end{eqnarray}%
\begin{eqnarray}
&&a_{3}\delta_{j-1}^{n+1}-a_{3}\delta_{j+1}^{n+1}-a_{1}\phi _{j-1}^{n+1}-a_{2}\phi
_{j}^{n+1}-a_{1}\phi _{j+1}^{n+1}  \label{k2} \\
&=&-a_{3}\delta_{j-1}^{n}+a_{3}\delta_{j+1}^{n}+a_{1}\phi _{j-1}^{n}+a_{2}\phi
_{j}^{n}+a_{1}\phi _{j+1}^{n}  \notag
\end{eqnarray}%
gives
\begin{eqnarray*}
&&\xi ^{n+1}\left[ A_1\left( 2a_{1}\cos \varphi +a_{2}\right) +\frac{A_2k\mu_3 }{2}%
\left( 2a_{4}\cos \varphi +a_{5}\right) -i\lambda k\varepsilon A_1a_{3}\sin
\varphi \right]  \\
&=&\xi ^{n}\left[ A_1\left( 2a_{1}\cos \varphi +a_{2}\right) -\frac{A_2k\mu_3 }{2}%
\left( 2a_{4}\cos \varphi +a_{5}\right) -i\lambda k\varepsilon A_1a_{3}\sin
\varphi \right] 
\end{eqnarray*}%
\begin{equation*}
\frac{\xi ^{n+1}}{\xi ^{n}}=\frac{\left[ A_1\left( 2a_{1}\cos \varphi
+a_{2}\right) -\frac{A_2k\mu_3 }{2}\left( 2a_{4}\cos \varphi +a_{5}\right)
-i\lambda k\varepsilon A_1a_{3}\sin \varphi \right] }{\left[ A_1\left(
2a_{1}\cos \varphi +a_{2}\right) +\frac{A_2k\mu_3 }{2}\left( 2a_{4}\cos \varphi
+a_{5}\right) -i\lambda k\varepsilon A_1a_{3}\sin \varphi \right] }
\end{equation*}%
\begin{equation}
\rho=\frac{\xi ^{n+1}}{\xi ^{n}}=\frac{X_{1}+iY}{X_{2}-iY}  \label{k3}
\end{equation}%
where%
\begin{eqnarray*}
X_{1} &=&A_1\left( 2a_{1}\cos \varphi +a_{2}\right) -\frac{A_2k\mu_3 }{2}\left(
2a_{4}\cos \varphi +a_{5}\right)  \\
X_{2} &=&A_1\left( 2a_{1}\cos \varphi +a_{2}\right) +\frac{A_2k\mu_3 }{2}\left(
2a_{4}\cos \varphi +a_{5}\right)  \\
Y &=&i\lambda k\varepsilon A_1a_{3}\sin \varphi 
\end{eqnarray*}%
and
\begin{eqnarray*}
&&\xi ^{n+1}\left[ -A_2\left( 2a_{1}\cos \varphi +a_{2}\right) -2iA_1a_{3}\sin
\varphi \right]  \\
&=&\xi ^{n}\left[ A_2\left( 2a_{1}\cos \varphi +a_{2}\right) +2iA_1a_{3}\sin
\varphi \right] 
\end{eqnarray*}%
\begin{equation*}
\frac{\xi ^{n+1}}{\xi ^{n}}=\frac{A_2\left( 2a_{1}\cos \varphi +a_{2}\right)
+2iA_1a_{3}\sin \varphi }{-A_2\left( 2a_{1}\cos \varphi +a_{2}\right)
-2iA_1a_{3}\sin \varphi }
\end{equation*}%
\begin{equation}
\rho =\frac{\xi ^{n+1}}{\xi ^{n}}=\frac{X_{3}+iZ}{X_{4}-iZ}  \label{k4}
\end{equation}%
\begin{eqnarray*}
X_{3} &=&A_2\left( 2a_{1}\cos \varphi +a_{2}\right)  \\
X_{4} &=&-A_2\left( 2a_{1}\cos \varphi +a_{2}\right)  \\
Z &=&2iA_1a_{3}\sin \varphi 
\end{eqnarray*}
It can be concluded from both (\ref{k3}) and (\ref{k4}) that $\left \vert \rho \right \vert$ is less than or equal to $1$. Thus, the proposed method method is unconditionally stable.

\section{Numerical Illustrations}
The numerical solutions for some IBVPs set up with the Gardner equation are summarized in this section. The accuracy of the results determined by the extended B-spline collocation method is discussed by examining graphical representations, measuring the error between the numerical and the analytical solutions and the preservation of conservation laws. The error of the numerical solution is measured by using the discrete maximum norm defined as 

\begin{equation*}
L_{\infty }(t)=\left \vert u(x_j,t)-U(x_j,t)\right \vert _{\infty }=\max \limits_{j}\left
\vert u(x_j,t)-U_{j}^{n}\right \vert
\end{equation*}%
where $U_j^n$ and $u(x_j,t)$ are numerical and analytical solutions at the discrete time $t$. 

The conservation laws can also be indicators of the validity of the proposed algorithms even when the analytical solution does not exist. The conservation laws of the Gardner equation  
\begin{equation}
\begin{aligned}
M&=\int\limits_{-\infty}^{\infty}{udx} \\
E&=\int\limits_{-\infty}^{\infty}{u^2dx} \\
H&=\int\limits_{-\infty}^{\infty}{\dfrac{\mu_1 u^3}{3}+\dfrac{\mu_2 u^4}{6}-\mu_3 (u_x)^2dx}
\end{aligned}
\end{equation} 
are expected to keep their initial quantities during numerical simulations \cite{hamdi1}. The relative changes of these quantities at a discrete time $t>0$ are measured by using $C(M_t)$, $C(E_t)$ and $C(H_t)$ defined as
\begin{equation}
\begin{aligned}
C(M_t)&=\left | \frac{M_t-M_0}{M_0} \right | \\
C(E_t)&=\left | \frac{E_t-E_0}{E_0} \right | \\
C(H_t)&=\left| \frac{H_t-H_0}{H_0} \right|
\end{aligned}
\end{equation}
where $M_t$, $E_t$, $H_t$, ($t\geq 0$) are the measured quantities at the time $t$.

\subsection{Propagation of Initial Single Positive Pulse}
\begin{figure}[h]
	\centering
		\includegraphics{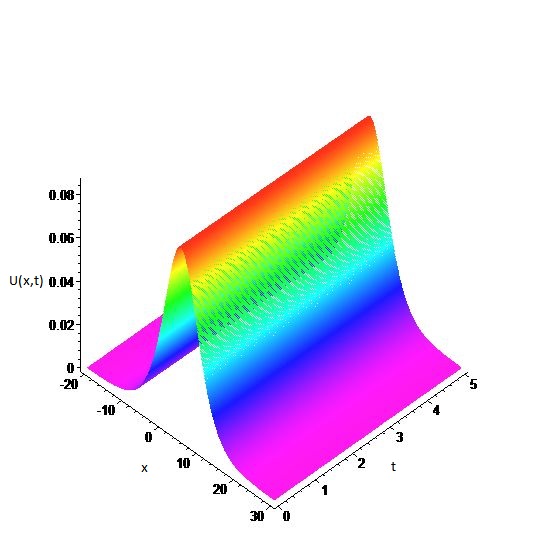}
	\caption{Propagation of a single solitary wave}	
	\label{fig:sech}
\end{figure}

In the first numerical illustration, we study propagation of an initial single pulse with positive amplitude. The equation parameters are chosen as $\mu _{1}=4$, $\mu _{2}=-3$ and $\mu _{3}=1$. The initial data are determined from the exact solution \cite{wazwaz1} 
\begin{equation*}
u(x,t)=\frac{2}{12+3\sqrt{14}\cosh (-\frac{x}{3}+\frac{5}{3}+\frac{t}{27})}
\end{equation*}
by assuming $t=0$. Since the exact solution approaches zero as $x$ approaches infinity or minus infinity, the choice of homogeneous boundary conditions is compatible with the solution. The artificial interval $[-20,30]$ is chosen for the numerical simulation and the designed algorithm is run up the simulation ending time $t=5$ with the fixed time discretization parameter $\Delta t=0.1$ and various spatial discretization numbers $N$. The propagation of a positive single solitary wave is depicted in Fig \ref{fig:sech}.

\begin{figure}[h]
    \subfigure[Error distribution for $\lambda =0$ at the simulation terminating time]{
   \includegraphics[scale =0.5] {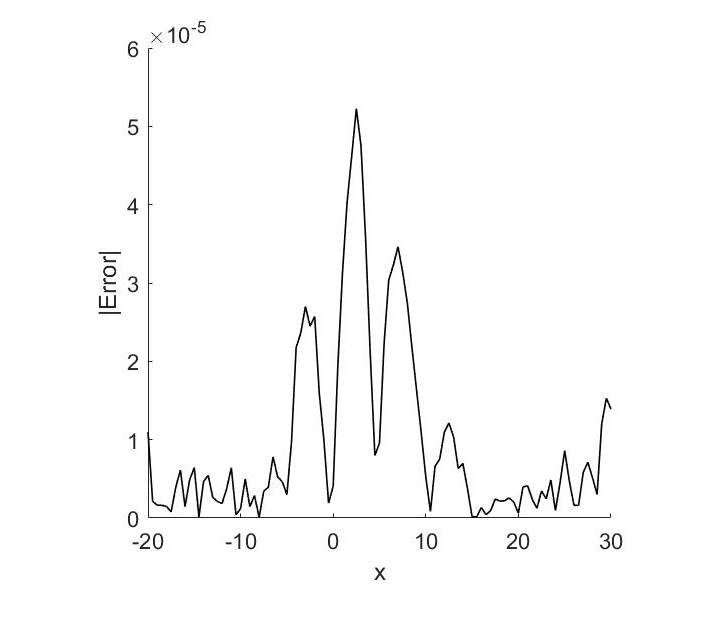}
   \label{fig:1a}
 }
 \subfigure[Error distribution for $\lambda=-0.00840$ at the simulation terminating time]{
   \includegraphics[scale =0.7] {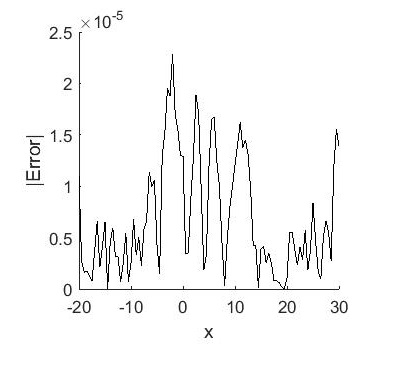}
   \label{fig:1b}
 }
 \caption{Error of the numerical error at the time $t=5$}
\end{figure}

The error distribution for the extension parameter $\lambda =0$ and the discretization parameters $\Delta t=0.1$ and $N=100$ at the simulation ending time $t=5$ is depicted in Fig \ref{fig:1a}. Usage of the same discretization parameters for optimum extension parameter $\lambda = -0.00840$ results in the error distribution at the simulation ending time as given in Fig \ref{fig:1b}. A simple comparison shows that the results are improved approximately two times when the optimum extension parameter is used in the algorithm. 
The discretized maximum error norms for various values of grid size are summarized in Table \ref{ex1norm}. In the meanwhile, the algorithm seeks for the optimal value of the extension parameter $\lambda$ by comparing the discrete maximum error norm at each $\lambda$. Even though we do not observe an improvement in decimal digits in the results, using the optimum extension parameter $\lambda$ generates improved results for all choices of $N$. 

\begin{table}[h]
\caption{Error norms for various choices of extension parameter}
	\centering
		\begin{tabular}{lllll}
		\hline \hline
		$N$ & $L_{\infty }(2.5)(\lambda =0)$ & $L_{\infty }(2.5)($various $\lambda)$ & $%
L_{\infty }(5)(\lambda =0)$ & $L_{\infty }(5)($various $\lambda)$ \\
$100$ & $3.2726\times 10^{-5}$ & $(\lambda=-0.00840)1.2330\times 10^{-5}$ & $%
5.22606\times 10^{-5}$ & $2.2789\times 10^{-5}$ \\ 
$200$ & $2.0537\times 10^{-5}$ & $(\lambda=-0.00280)1.4819\times 10^{-5}$ & $%
1.91604\times 10^{-5}$ & $1.9119\times 10^{-5}$ \\ 
$300$ & $1.4428\times 10^{-5}$ &$(\lambda=-0.00094)1.2509\times 10^{-5}$  & $1.70403\times 10^{-5}$ &  $1.6944\times 10^{-5}$\\ 
$400$ & $1.4452\times 10^{-5}$ & $(\lambda=-0.00178)1.4440\times 10^{-5}$ & $1.61150\times 10^{-5}$ & $1.5872\times 10^{-5}$ \\ \hline	\hline	
		\end{tabular}
		\label{ex1norm}
\end{table}
The conservation laws are required to preserve their initial values as time proceeds during the simulation. The initial values of these laws are calculated by using Maple, Table \ref{ex1cl}. The absolute relative changes of conservation laws are obtained at least six decimal digits at the simulation ending time $t=5$. These preservation rates can be accepted as indicators of a valid algorithm. 
\begin{table}[h]
\caption{Calculated conservation laws and their absolute relative changes}
	\centering
		\begin{tabular}{lllllll}
		\hline \hline
	$N$ & $M_{0}$ & $E_{0}$ & $H_{0}$ & $C(M_{5})$ & $C(E_{5})$ & $C(H_{5})
$ \\ 
$100$ & $1.0445$ & $0.0601$ & $0.0040$ & $5.4748\times 10^{-6}$ & $%
3.8176\times 10^{-8}$ & $1.5233\times 10^{-6}$ \\ 
$200$ & $1.0445$ & $0.0601$ & $0.0040$ & $3.2669\times 10^{-6}$ & $%
5.1126\times 10^{-8}$ & $1.7003\times 10^{-6}$ \\ 
$300$ & $1.0445$ & $0.0601$ & $0.0040$ & $2.4190\times 10^{-7}$ & $%
2.1767\times 10^{-8}$ & $2.8351\times 10^{-6}$ \\ 
$400$ & $1.0445$ & $0.0601$ & $0.0040$ & $1.3753\times 10^{-6}$ & $%
2.0910\times 10^{-10}$ & $3.3939\times 10^{-6}$ \\
	 \hline	\hline	
		\end{tabular}
		\label{ex1cl}
\end{table}
\subsection{Propagation of Kink-like Wave}
\begin{figure}[h]
	\centering
		\includegraphics{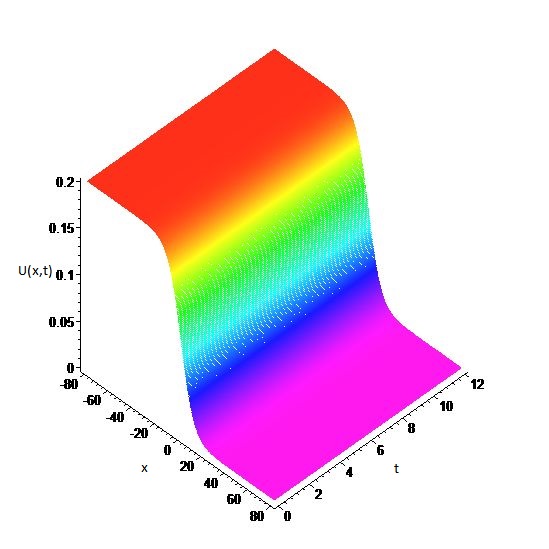}
	\caption{Propagation of a single solitary wave}	
	\label{fig:tanh}
\end{figure}

Kink-like wave solution of the Gardner equation is 
\begin{equation}
u(x,t)=\dfrac{1}{10}-\dfrac{1}{10}\tanh{(\dfrac{\sqrt{30}}{60}(x-\dfrac{1}{30}t))} \label{tanh}
\end{equation}
for the equation parameters $\mu_1=1$, $\mu_2=-5$ and $\mu_3=1$\cite{wazwazbook}. This wave propagates along the $x-$axis with the speed $1/30$. The initial data are generated from the analytical solution by assuming $t=0$. Since the analytical solution disappears as $x\rightarrow \infty $ and $u(x,t)\rightarrow 0.2$ as $x\rightarrow - \infty$, the homogeneous Neumann data are compatible. The artificial problem interval is chosen as $[-80,80]$ and the designed algorithm is run up to the ending time $t=12$ with $\Delta t =0.1$ and various grid numbers used in the space domain. The plot Fig \ref{fig:tanh} is a summary of the propagation in this finite interval.  

The designed algorithm is scanned the extension parameter in $[-1,1]$ with step size $\Delta \lambda =0.000001$ for an optimum value to improve the results. Optimum extension parameters are determined for all choices of the number of grids, Table \ref{ex2norm}. Determination of the optimum extension parameter reduces the maximum absolute error to the half for all choices of $N$.

\begin{table}[h]
\caption{Error norms for optimum choices of extension parameter}
	\centering
		\begin{tabular}{lllll}
		\hline \hline
		$N$ & $L_{\infty }(4)(\lambda =0)$ & $L_{\infty }(4)($various $\lambda)$ & $%
L_{\infty }(12)(\lambda =0)$ & $L_{\infty }(12)($various $\lambda)$ \\
$100$ & $8.4150\times 10^{-6}$ & $(\lambda=-0.01850)3.8974\times 10^{-6}$ & $%
2.3158\times 10^{-5}$ & $(\lambda=-0.01850)1.2330\times 10^{-5}$ \\ 
$200$ & $2.1207\times 10^{-6}$ & $(\lambda=-0.00574)1.0194\times 10^{-6}$ &  $5.9956\times
10^{-6}$ & $(\lambda=-0.00574)2.9662\times 10^{-6}$ \\ 
$400$ & $5.3296\times 10^{-7}$ & $(\lambda=-0.00115)2.5440\times 10^{-7}$ &  $1.5016\times
10^{-6}$ & $(\lambda=-0.00115)7.7413\times 10^{-7}$ \\ 
$600$ & $2.2377\times 10^{-7}$ & $(\lambda=-0.00057)1.1335\times 10^{-7}$ &$6.6655\times
10^{-7} $ & $(\lambda=-0.00057)3.3921\times 10^{-7}$ \\ 
$800$ & $1.4601\times 10^{-7}$ & $(\lambda=-0.00024)6.3749\times 10^{-8}$ &  $5.2835\times 10^{-6}$ & 
$(\lambda=-0.00024)5.2779\times 10^{-6}$ \\
 \hline	\hline	
		\end{tabular}
		\label{ex2norm}
\end{table}

The initial values of the conservation laws are calculated by integration the quantities by substituting the initial data of the IBVP. The absolute relative changes of the conservation laws indicate a reliable solution in Table \ref{ex2cl}.
\begin{table}[h]
\caption{Calculated conservation laws and their absolute relative changes}
	\centering
		\begin{tabular}{lllllll}
		\hline \hline
$100$ & $16.1599$ & $3.0129$ & $0.0979$ & $4.9504\times 10^{-3}$ & $%
5.3104\times 10^{-3}$ & $5.4423\times 10^{-3}$ \\ 
$200$ & $16.0799$ & $2.9969$ & $0.0974$ & $4.9751\times 10^{-3}$ & $%
5.3388\times 10^{-3}$ & $5.4721\times 10^{-3}$ \\ 
$400$ & $16.0399$ & $2.9889$ & $0.0972$ & $4.9875\times 10^{-3}$ & $%
5.3531\times 10^{-3}$ & $5.4871\times 10^{-3}$ \\ 
$600$ & $16.0266$ & $2.9862$ & $0.0971$ & $4.9916\times 10^{-3}$ & $%
5.3578\times 10^{-3}$ & $5.4922\times 10^{-3}$ \\ 
$800$ & $16.0199$ & $2.9849$ & $0.0970$ & $4.9938\times 10^{-3}$ & $%
5.3603\times 10^{-3}$ & $4.9481\times 10^{-3}$ \\
	 \hline	\hline	
		\end{tabular}
		\label{ex2cl}
\end{table}
\begin{figure}[h]
    \subfigure[Error distribution for $\lambda =0$ at the simulation terminating time]{
   \includegraphics[scale =0.6] {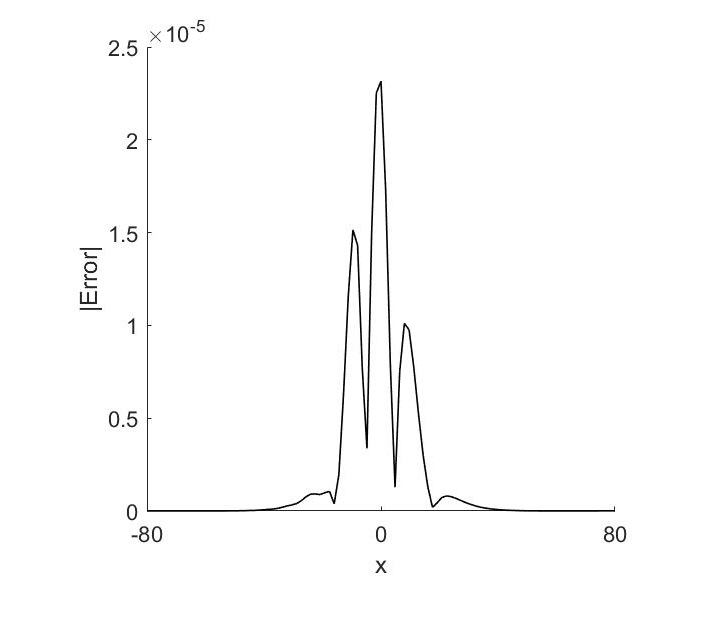}
   \label{fig:2a}
 }
 \subfigure[Error distribution for $\lambda =-0.01850$ at the simulation terminating time]{
   \includegraphics[scale =0.7] {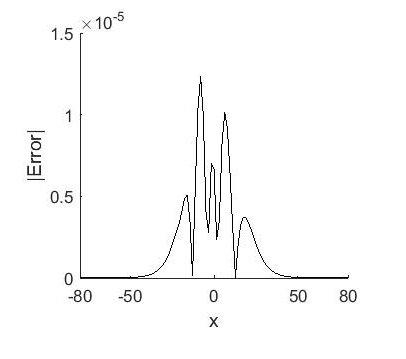}
   \label{fig:2b}
 }
 \caption{Error of the numerical error at the time $t=12$}
\end{figure}
\subsection{Wave Generation from an Initial Pulse}
\begin{figure}[hf]
    \subfigure[Initial data]{
   \includegraphics[scale =0.65] {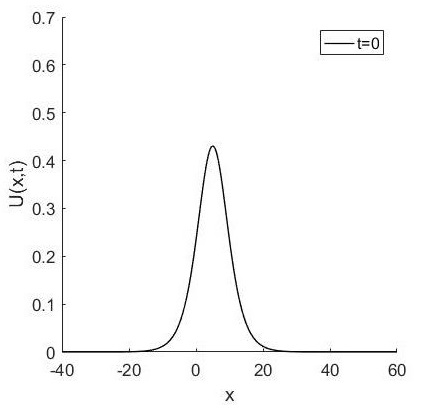}
   \label{fig:3a}
 }
   \subfigure[$t=5$]{
   \includegraphics[scale =0.65] {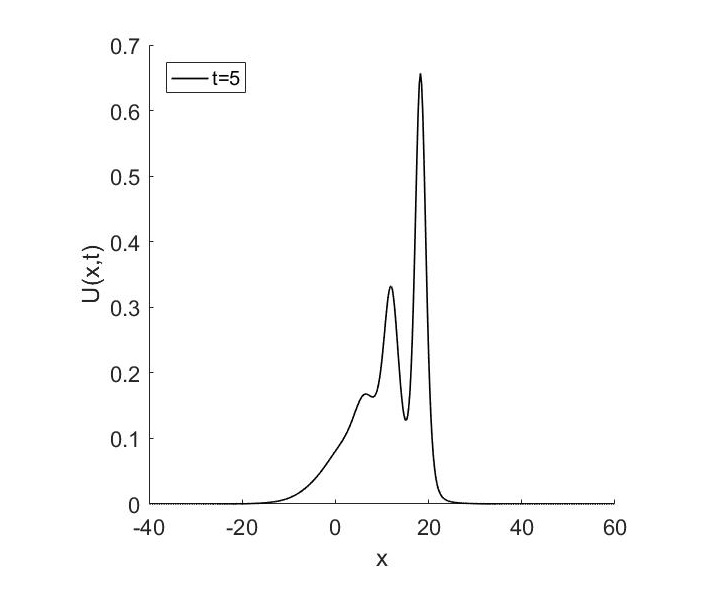}
   \label{fig:3b}
 }
 \subfigure[$t=10$]{
   \includegraphics[scale =0.65] {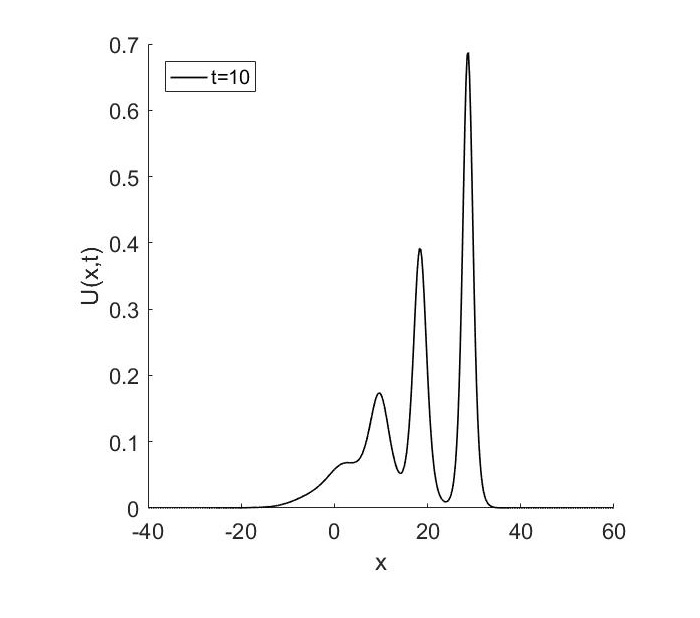}
   \label{fig:3c}
 }
 \subfigure[$t=15$]{
   \includegraphics[scale =0.65] {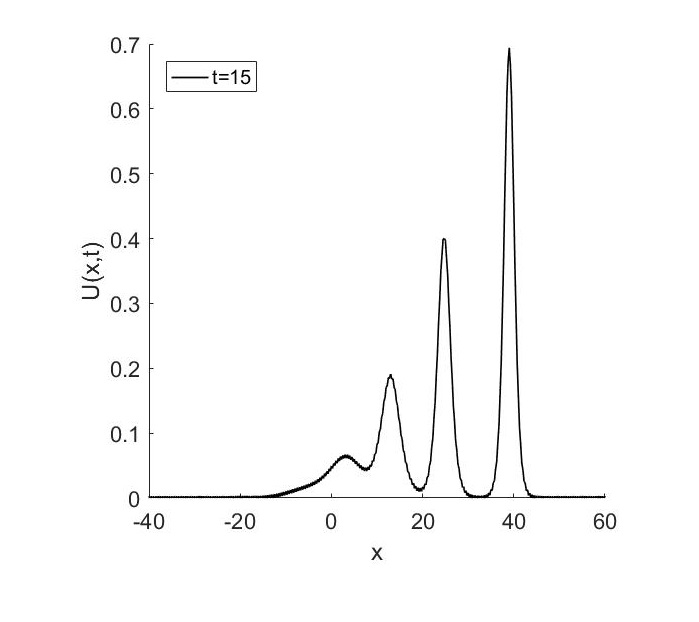}
   \label{fig:3d}
 }
 \caption{Wave generation from an initial positive pulse}
\end{figure}

\noindent
The perturbed Gardner equation of the form
\begin{equation}
u_{t}+\mu_1 uu_{x}+\mu_2 u^2u_{x}+\mu_3 u_{xxx}+\epsilon=0 \label{gardnerp}
\end{equation}
for some nonzero real $\epsilon$ can be useful to study the wave generation from an initial positive pulse. The initial data are chosen as
\begin{equation}
u(x,t)=\dfrac{2}{3}\dfrac{1}{4+\sqrt{14}\cosh{\left(\dfrac{x}{3}-\dfrac{5}{3}\right)}}
\end{equation}
with the parameter choice $\mu_1 =10$, $\mu_2 =-3$ and $\mu_3 =1$ in the Gardner Equation (\ref{gardnerp}). The designed program is run in the finite interval $[-40,60]$ with the parameters $N=400$ and $\Delta t=0.01$ up to the $t=15$. The initial positive pulse of height $0.4305$ positioned at $x=5$ propagates along the $x$-axis, Fig \ref{fig:3a}. When the propagation time reaches $t=5$, a frontier wave of height $0.6568$ is positioned at $x=18.25$, Fig \ref{fig:3b}. The first follower solitary wave is clearly observable at this time. This follower is of height $0.3318$ and it is positioned at $x=12$. Even though it is not clearly observable, a bulge at the left of the first follower can be evaluated as an indicator of a second follower wave. The first three waves are clearly observable at $t=10$, Fig \ref{fig:3c}. The height of the frontier reaches $0.6871$ and its peak point is positioned at $x=28.75$. The peak of the first follower wave of height $0.3913$ is positioned at $x=18.25$. At the left of the first follower, almost completely formed second follower of height $0.1736$ is positioned at $x=9.75$. The height of the frontier is measured as $0.6941$ at the time $t=15$, Fig \ref{fig:3d}. Its peak is positioned at $x=39$ at this time. The height of the first follower wave reaches $0.3998$ and its peak position is measured as $x=24.75$.  The second follower of height $0.1910$ is positioned at $x=13$. The bulge appearing at the left of the second follower is the indicator of formation of one more solitary.        

Even though we can not compare the results with the analytical solutions, the absolute relative changes of the conservation laws is a good indicator to evaluate the efficiency of the proposed method. The initial values of the conservation laws and their absolute relative changes are summarized in Table \ref{ex3cl}. The absolute relative change of the first conservation law is only in six decimal digits, as the second one in four decimal digits and the third one in three decimal digits at the simulation ending time $t=15$.

\begin{table}[hf]
\caption{Calculated conservation laws and their absolute relative changes}
	\centering
		\begin{tabular}{lllllll}
		\hline \hline
		$t$ & $M_0$ & $E_0$ & $H_0$ & $C(M_15)$ & $C(E_15)$ & $C(H_15)$ \\
$5$ & $5.2255$ & $1.5033$ & $1.5994$ & $8.0719\times 10^{-7}$ & $%
3.0588\times 10^{-5}$ & $1.2886\times 10^{-3}$ \\ 
$10$ & $5.2255$ & $1.5033$ & $1.5994$ & $2.7652\times 10^{-6}$ & $%
4.1342\times 10^{-5}$ & $1.8485\times 10^{-3}$ \\ 
$15$ & $5.2255$ & $1.5033$ & $1.5994$ & $7.0380\times 10^{-6}$ & $%
6.1132\times 10^{-4}$ & $2.1571\times 10^{-3}$\\
	 \hline	\hline	
		\end{tabular}
		\label{ex3cl}
\end{table}

\section{Conclusion}
For the sake of the improve the results obtained by the cubic B-splines, the extended cubic B-spline based collocation method is constructed for some initial boundary value problems for the Gardner Equation. The scan of the optimum extension parameter in the interval $[-1,1]$ gives opportunity to obtain more reliable results when compare classical cubic B-splines. Since the extended cubic b-splines have only first and second order derivatives, the reduction of the order of the third order derivative is required. Thus, the coupled system of nonlinear PDEs is obtained. The extended B-spline function are used to approximate the solutions of this system. Following the spatial discretization, the linearization procedure is followed. At the end, the time integration is done Crank-Nicolson method. Von-Neumann stability analysis shows that the suggested algorithm is unconditionally stable.

The first two examples give opportunity to measure the error between the analytical and the numerical solutions by calculating maximum error norms for various choices of the discretization parameters. Both graphical representations and the absolute relative changes of the conservation laws are indicators of a reliable and valid method. The results are improved by determining the optimum value of the extension parameter. 

In the third example, a non analytical problem simulating wave generation from an initial single solitary wave is studied. The proposed algorithm simulates the expected results successfully. The absolute relative changes of the conservation laws confirm the valid results.

\section{Acknowledgements}
\textsl{This study is a part of the project with number 2016/19052 supported by Eskisehir Osmangazi University Scientific Research Projects Committee and was partially presented at 3rd International Conference on Pure and Applied Sciences, Dubai, 2017.}

\end{document}